\documentclass[12pt]{amsart}
\usepackage{amssymb,amsmath,amsthm}
\usepackage[english]{babel}
\usepackage{epsfig}
\usepackage{pstricks}
\usepackage{comment}

\numberwithin{equation}{section}

\newcommand{\Z}{ \mathbb Z}
\newcommand{\N}{ \mathbb N}

\newcommand{\C}{ \mathbb C}
\newcommand{\p}[1]{\mathcal{#1}}
\newcommand{\s}[1]{\mathbf{#1}}
\newcommand{\vp}{\varphi}
\newcommand{\set}[1]{\{#1\}}
\newcommand{\detq}{{\textstyle \det_q}}
\newcommand{\detqq}{{\textstyle \det_{\s q}}}
\DeclareMathOperator{\inv}{inv}
\DeclareMathOperator{\comm}{{comm}}
\DeclareMathOperator{\cf}{{cf}}
\DeclareMathOperator{\rtq}{{rq}}

\newtheoremstyle{thm}
  {9pt}{9pt}{\itshape}{}{\bfseries}{}{.5em}{}
\theoremstyle{thm}
\newtheorem{thm}{Theorem}[section]

\newtheorem{lemma}[thm]{Lemma}
\newtheorem{prop}[thm]{Proposition}

\newtheoremstyle{defin}
  {9pt}{9pt}{}{}{\bfseries}{}{.5em}{}
\theoremstyle{defin}

\newtheoremstyle{exm}
  {9pt}{9pt}{}{}{\scshape}{}{.5em}{}
\theoremstyle{exm}
\newtheorem{exm}[thm]{Example}
\newtheorem{rmk}[thm]{Remark}
\newtheoremstyle{prf}
  {}{}{}{}{\itshape}{:}{.5em}{}
\theoremstyle{prf}
\newtheorem*{prf}{Proof}

\title{Non-commutative Sylvester's determinantal identity}
\author{Matja\v z Konvalinka}
\date{\today}

\begin{document}

\begin{abstract}
 Sylvester's identity is a classical determinantal identity with a straightforward linear algebra proof. We present a new, combinatorial proof of the identity, prove several non-commutative versions, and find a $\beta$-extension that is both a generalization of Sylvester's identity and the $\beta$-extension of the MacMahon master theorem.
\end{abstract}

\maketitle

\section{Introduction} \label{intro}

\subsection{Classical Sylvester's determinantal identity.} Sylvester's identity is a classical determinantal identity that is usually written in the form used by Bareiss (\cite{bareiss}).

\begin{thm}[Sylvester's identity] \label{intro2}
 Let $A$ denote a matrix $(a_{ij})_{m \times m}$; take $n < i,j \leq m$ and define
 $$A_0 = \begin{pmatrix} a_{11} & a_{12} & \cdots & a_{1n} \\
			   a_{21} & a_{22} & \cdots & a_{2n} \\
			   \vdots & \vdots & \ddots & \vdots \\
			   a_{n1} & a_{n2} & \cdots & a_{nn} \end{pmatrix}, \quad 
 a_{i*} = \begin{pmatrix} a_{i1} & a_{i2} & \cdots & a_{in} \end{pmatrix}, \quad a_{*j} = \begin{pmatrix} a_{1j} \\ a_{2j} \\ \vdots \\ a_{nj} \end{pmatrix},$$
 $$b_{ij} = \det \begin{pmatrix} A_0 & a_{*j} \\ a_{i*} & a_{ij} \end{pmatrix}, \quad B = (b_{ij})_{n+1 \leq i,j \leq m}$$
 Then
 $$\det A \cdot (\det A_0)^{m-n-1} = \det B. \eqno \qed$$
\end{thm}

\begin{exm}
 If we take $n = 1$ and $m = 3$, the Sylvester's identity says that
 $$(a_{11}a_{22}a_{33}-a_{11}a_{32}a_{23}-a_{21}a_{12}a_{33}+a_{21}a_{32}a_{13}+a_{31}a_{12}a_{23}-a_{31}a_{22}a_{13})a_{11} =$$
 $$ = \begin{vmatrix} a_{11}a_{22} - a_{21} a_{12} & a_{11}a_{23} - a_{21} a_{13} \\ a_{11}a_{32} - a_{31} a_{12} & a_{11}a_{33} - a_{31} a_{13} \end{vmatrix}.$$
\end{exm}

The Sylvester's identity has been intensely studied, mostly in the algebraic rather than combinatorial context. The crucial step was made by Krob and Leclerc~\cite{krob}, who found a quantum version. Since then, Molev found several far-reaching extensions to Yangians, including other root systems~\cite{molev1,molev2} (see also~\cite{hm}).
			   
\subsection{Main result.} In this paper, we will find a new combinatorial proof of the classical Sylvester's identity and find a multiparameter right-quantum analogue. We use the techniques developed in \cite{konvalinka}.

\medskip

Fix non-zero complex numbers $q_{ij}$ for $1 \leq i < j \leq m$. We call a matrix $A$ $\s q$-right-quantum if
\begin{eqnarray}
a_{jk}  a_{ik} &  = &  q_{ij}   a_{ik}  a_{jk}  \ \
\text{for all} \ \ i < j,\\
a_{ik}  a_{jl}   -   q_{ij}^{-1}   a_{jk}  a_{il} &  = & 
q_{kl}q_{ij}^{-1}   a_{jl}  a_{ik}   -  
q_{kl}   a_{il}  a_{jk} \ \
\text{for all} \ \ i < j, \ k <l.
\end{eqnarray}

In the next section, we will the define the concept of a $\s q$-determinant of a square matrix. We will have

$$\detqq(I-A)= \sum_{J \subseteq [m]} (-1)^{|J|} \detqq A_J,$$
where
$$\detqq A_J = \sum_{\sigma \in S_J} \left( \prod_{p<r \colon j_p>j_r} q_{j_r j_p}^{-1}\right)a_{\sigma(j_1)j_1} \cdots a_{\sigma(j_k)j_k}$$
for $J=\set{j_1<j_2<\ldots<j_k}$.

\medskip

Our main theorem is the following.

\begin{thm}[$\s q$-right-quantum Sylvester's determinant identity] \label{intro3}
 Let $A=(a_{ij})_{m \times m}$ be a $\s q$-right-quantum matrix, and choose $n < m$. Let $A_0,a_{i*},a_{*j}$ be defined as above, and let
 $$c^{\s q}_{ij} = - \detqq\!\!\!{}^{-1}(I-A_0) \cdot \detqq \begin{pmatrix} I-A_0 & -a_{*j} \\ -a_{i*} & -a_{ij} \end{pmatrix}, \quad C^{\s q}=(c_{ij}^{\s q})_{n+1 \leq i,j \leq m}.$$
 Suppose $q_{ij}=q_{i'j'}$ for all $i,i' \leq n$ and $j,j' > n$. Then
 $$\detqq\!\!\!{}^{-1}(I-A_0) \cdot \detqq(I-A)=\detqq (I-C^{\s q}).$$
\end{thm}

The determinant $\detqq(I-A_0)$ does not commute with other determinants in the definition of $c_{ij}^{\s q}$, so the identity cannot be written in a form analogous to Theorem \ref{intro2}. See Remark \ref{cfqij11} for a discussion of the necessity of the condition $q_{ij}=q_{i'j'}$ for $i,i' \leq n$, $j,j' > n$.

\medskip

The proof roughly follows the pattern of the proof of the main theorem in \cite{konvalinka}. First we show a new, combinatorial proof of the classical Sylvester's identity (Sections \ref{sylv} and \ref{comm}). Then we adapt the proof to simple non-commutative cases -- the Cartier-Foata case (Section \ref{cf}) and the right-quantum case (Section \ref{rq}). We extend the results to cases with a weight (Sections \ref{cfq} and \ref{qrq}) and to multiparameter weighted cases (Sections \ref{cfqij} and \ref{qrqij}). We also present a $\beta$-extension of Sylvester's identity in Section \ref{beta}.

\section{Algebraic framework} \label{alg}

\subsection{Words and matrices.} We will work in the $\C$-algebra $\p A$ of formal power series in non-commuting variables $a_{ij}$, $1 \leq i,j \leq m$. Elements of $\p A$ are infinite linear combinations of words
in variables $a_{ij}$ (with coefficients in~$\C$). In most cases we will take elements of $\p A$ modulo some ideal~$\p I$ generated by a finite number of quadratic relations. For example, if~$\p I_{\comm}$ is generated by $a_{ij}a_{kl} = a_{kl}a_{ij}$ for all $i,j,k,l$, then $\p A/\p I_{\comm}$ is the symmetric algebra (the free commutative algebra with variables $a_{ij}$).

\medskip

We abbreviate the product $a_{\lambda_1\mu_1}\cdots a_{\lambda_\ell\mu_\ell}$
to $a_{\lambda,\mu}$ for $\lambda=\lambda_1\cdots \lambda_\ell$ and
$\mu=\mu_1\cdots \mu_\ell$, where~$\lambda$ and~$\mu$ are regarded as words
in the alphabet $\set{1,\ldots,m}$. For such a word $\nu=\nu_1 \cdots \nu_\ell$,
define the \emph{set of inversions}
$$\p I(\nu) \, = \, \set{(i,j) \colon i < j, \nu_i > \nu_j},$$
and let~$\inv\nu = |\p I(\nu)|$ be the \emph{number of inversions}.

\subsection{Determinants.} Let $B=(b_{ij})_{n \times n}$ be a square matrix with entries
in $\p A$, i.e.\hspace{-0.07cm} $b_{ij}$'s are linear combinations of words in $\p A$. To define the determinant of $B$, expand the terms of
$$\sum_{\sigma \in S_n} (-1)^{\inv(\sigma)} b_{\sigma_1{1}} \cdots b_{\sigma_n n},$$
and weight a word $a_{\lambda,\mu}$ with a certain weight $w(\lambda,\mu)$. The resulting expression will be called the \emph{determinant} of $B$ (with respect to $\p A$). In the usual commutative case, all weights are equal to $1$.

\medskip

In all cases we consider we have $w(\varnothing,\varnothing)=1$. Therefore
$$\frac1{\det(I-A)} \, = \, \frac1{1 - \Sigma} \, = \, 1 \, + \, \Sigma \,
+ \, \Sigma^2 \, + \, \ldots\,,
$$
where~$\Sigma$ is a certain finite sum of words in~$a_{ij}$ and both the left and the right inverse of $\det(I-A)$ are equal to the infinite sum on the right. We can use the fraction notation as above in non-commutative situations.

\subsection{Paths.} We will consider \emph{lattice steps} of the form $(x,i) \to (x+1,j)$ for some $x,i,j \in \Z$,
$1 \leq i,j \leq m$.  We think of $x$ being drawn along the $x$-axis, increasing from left to right, and refer to~$i$ and~$j$ as the \emph{starting height} and \emph{ending height}, respectively. We identitfy the step $(x,i) \to (x+1,j)$ with the variable $a_{ij}$. Similarly, we identify a finite sequence of steps with a word in the alphabet $\set{a_{ij}}$, $1 \leq i,j \leq m$, i.e.\hspace{-0.07cm} with an element of the algebra $\p A$. If each step in a sequence starts at the ending point of the previous step, we call such a sequence a \emph{lattice path}. A lattice path with starting height $i$ and ending height $j$ will be called a path from $i$ to $j$.

\begin{exm} \label{alg1}
 The following is a path from $4$ to $4$.
 \begin{figure}[ht!]
 \begin{center}
  \epsfig{file=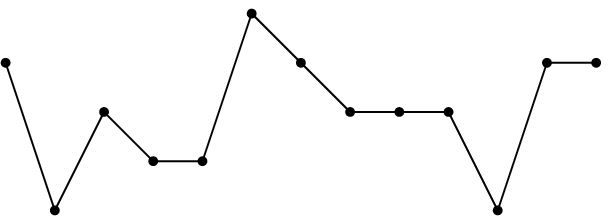,height=1.5cm}
  \caption{Representation of the word $a_{41}a_{13}a_{32}a_{22}a_{25}a_{54}a_{43}a_{33}a_{33}a_{31}a_{14}a_{44}$.}
  \label{fig2}
 \end{center}
 \end{figure}
\end{exm}

Recall that the $(i,j)$-th entry of $A^k$ is the sum of all paths of length $k$ from $i$ to $j$. Since
$$(I-A)^{-1} = I + A + A^2 + \ldots,$$
the $(i,j)$-th entry of $(I-A)^{-1}$ is the sum of all paths (of any length) from $i$ to $j$.

\section{Non-commutative Sylvester's identity} \label{sylv}

As in Section \ref{intro}, choose $n < m$, and denote the matrix $(a_{ij})_{m \times m}$ by $A$ and $(a_{ij})_{n \times n}$ by $A_0$.

\medskip

We will show a combinatorial proof of the non-commutative Sylvester's identity due to Gelfand and Retakh, see \cite{gelfand2}.

\begin{thm}[Gelfand-Retakh] \label{sylv2}
 Consider the matrix $C = (c_{ij})_{n+1 \leq i,j \leq m}$, where
 $$c_{ij}=a_{ij} + a_{i*} (I-A_0)^{-1} a_{*j}.$$
 Then
 $$(I-A)^{-1}_{ij} = (I-C)^{-1}_{ij}.$$
\end{thm}
\begin{prf}
 Take a lattice path $a_{ii_1}a_{i_1i_2}\cdots a_{i_{\ell-1}j}$ with $i,j>n$. Clearly it can be uniquely divided into paths $P_1, P_2, \ldots P_p$ with the following properties:
 \begin{itemize}
  \item the ending height of $P_i$ is the starting height of $P_{i+1}$
  \item the starting and the ending heights of all $P_i$ are strictly greater than $n$
  \item all intermediate heights are less than or equal to $n$
 \end{itemize}
 Next, note that
 $$c_{ij} = a_{ij} + a_{i*} (I-A_0)^{-1} a_{*j} = a_{ij} + \sum_{k,l \leq n} a_{ik} (I + A_0 + A_0^2 + \ldots)_{kl}a_{lj}$$
 is the sum over all non-trivial paths with starting height $i$, ending height $j$, and intermediate heights $\leq n$. This decomposition hence proves the theorem.\qed
\end{prf}

\begin{exm}
The following figure depicts the path from Example \ref{alg1} with a dotted line between heights $n$ and $n+1$, and the corresponding decomposition, for $n=3$.

 \begin{figure}[ht]
  \begin{center}
   \input{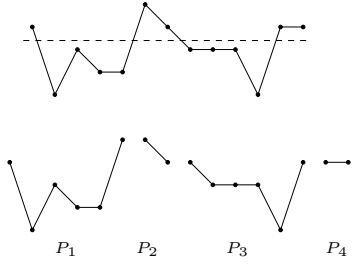}
  \end{center}
  \caption{The decomposition $(a_{41}a_{13}a_{32}a_{22}a_{25})(a_{54})(a_{43}a_{33}a_{33}a_{31}a_{14})(a_{44})$.}
  \label{fig1}
 \end{figure} 
\end{exm}

The theorem implies that
\begin{equation} \label{sylv1}
 (I-A)^{-1}_{n+1,n+1}(I-A^{n+1,n+1})^{-1}_{n+2,n+2}\cdots \left(I - \begin{pmatrix} A_0 & a_{*m} \\ a_{m*} & a_{mm} \end{pmatrix} \right)^{-1}_{mm} =
\end{equation}
$$=(I-C)^{-1}_{n+1,n+1}(I-C^{n+1,n+1})^{-1}_{n+2,n+2}\cdots (1-c_{mm})^{-1}.$$

Here $A^{n+1,n+1}$ is the matrix $A$ with the $(n+1)$-th row and column removed.

\medskip

In all the cases we will consider in the following sections, both the left-hand side and the right-hand side of this equation will be written in terms of determinants, as in the classical Sylvester's identity.

\section{Commutative case} \label{comm}

Recall that if $D$ is an invertible matrix with commuting entries, we have
$$\left(D^{-1}\right)_{ij} = (-1)^{i+j}\frac{\det D^{ji}}{\det D},$$
where $D^{ji}$ denotes the matrix $D$ without the $j$-th row and the $i$-th column. Apply this to \eqref{sylv1}: the numerators (except the last one on the left-hand side) and denominators (except the first one on both sides) cancel each other, and we get
\begin{equation} \label{comm1}
 \frac{\det(I-A_0)}{\det(I-A)} = \frac 1 {\det (I-C)}.
\end{equation}

\begin{prop} \label{comm2}
 For $i,j > n$ we have
 \begin{equation} \label{comm3}
  \delta_{ij} - c_{ij} = \frac{\det \begin{pmatrix} I - A_0 & -a_{*j} \\ -a_{i*} & \delta_{ij}-a_{ij} \end{pmatrix}}{\det(I-A_0)}.
 \end{equation}
\end{prop}
\begin{prf}
 Clearly we have
 $$(1-c_{ij})^{-1} = \left(\left( I - \begin{pmatrix} A_0 & a_{*j} \\ a_{i*} & a_{ij} \end{pmatrix}\right)^{-1}\right)_{ij},$$
 and by \eqref{comm1}, this is equal to
 $$\frac{\det(I-A_0)}{\det\left( I - \begin{pmatrix} A_0 & a_{*j} \\ a_{i*} & a_{ij} \end{pmatrix} \right)}.$$
 This finishes the proof for $i=j$, and for $i \neq j$ we have
 $$1-c_{ij} = \frac{\det \begin{pmatrix} I - A_0 & -a_{*j} \\ -a_{i*} & 1 - a_{ij} \end{pmatrix}}{\det(I-A_0)} = \frac{\det \begin{pmatrix} I - A_0 & -a_{*j} \\ -a_{i*} & - a_{ij} \end{pmatrix} + \det \begin{pmatrix} I - A_0 & 0 \\ -a_{i*} & 1 \end{pmatrix}}{\det(I-A_0)} = $$
 $$= \frac{\det \begin{pmatrix} I - A_0 & -a_{*j} \\ -a_{i*} & - a_{ij} \end{pmatrix} + \det(I-A_0)}{\det(I-A_0)} = \frac{\det \begin{pmatrix} I - A_0 & -a_{*j} \\ -a_{i*} & - a_{ij} \end{pmatrix}}{\det(I-A_0)} + 1. \eqno \qed$$ 
\end{prf}

\begin{prf}[of Theorem \ref{intro2}] 
 The proposition, together with \eqref{comm1}, implies that
 $$\frac{\det(I-A)}{\det(I-A_0)} = \det (I-C) = \det(I-A_0)^{n-m} \det B$$
 for
 $$b_{ij} = \det \begin{pmatrix} I - A_0 & -a_{*j} \\ -a_{i*} & \delta_{ij}-a_{ij} \end{pmatrix}, \quad B = (b_{ij})_{n+1 \leq i,j \leq m},$$
 which is Theorem \ref{intro2} for the matrix $I-A$.\qed
\end{prf}

\section{Cartier-Foata case} \label{cf}

A matrix $A$ is Cartier-Foata if
\begin{equation} \label{cf6}
 a_{ik}a_{jl}=a_{jl}a_{ik}
\end{equation}
for $i \neq j$, and right-quantum if
\begin{eqnarray}
 a_{jk}  a_{ik} & = &  a_{ik}a_{jk}  \ \
 \text{for all}  \ \ i \neq j ,
 \label{cf7} \\
 a_{ik}  a_{jl}  -  a_{jk}  a_{il}   & = &  
 a_{jl}  a_{ik}  - a_{il}  a_{jk} 
    \ \ \text{for all}  \ \ i \neq j , k \neq l . \label{cf8}
\end{eqnarray}

A Cartier-Foata matrix is also right-quantum, but the proofs tend to be much simpler for Cartier-Foata matrices.

\medskip

Note also that the classical definition of the determinant
$$\det B = \sum_{\sigma \in S_m} (-1)^{\inv \sigma}
 b_{\sigma_1 1}\cdots b_{\sigma_m m}$$
makes sense for a matrix $B = (b_{ij})_{m \times m}$ with entries generated by $a_{ij}$; in the language of Section \ref{alg}, we have $w(\lambda,\mu)=1$ for all words $\lambda,\mu$.

\medskip

Recall the following result (see Section \ref{final}).

\begin{prop} \label{cf1}
 If $A=(a_{ij})_{m \times m}$ is a Cartier-Foata matrix or a right-quantum matrix, we have
 $$\left(\frac1{I-A}\right)_{ij} \, = \, (-1)^{i+j} \frac1{\det(I-A)} \, \cdot \,
   \det\left(I-A\right)^{ji}$$
 for all $i,j$.\qed
\end{prop}

\begin{figure}[ht!]
 \begin{center}
  \epsfig{file=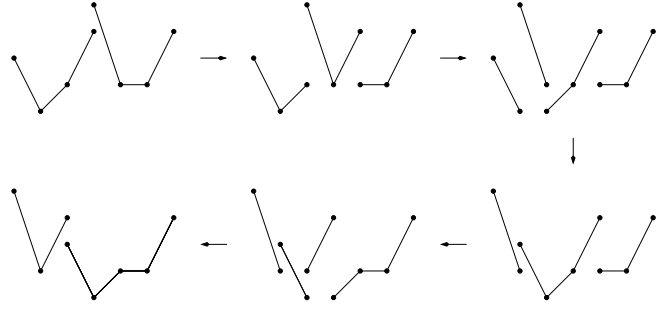,height=4cm}
  \caption{Transforming $a_{31}a_{12}a_{24}a_{52}a_{22}a_{24}$ into $a_{52}a_{24}a_{31}a_{12}a_{22}a_{24}$.}
  \label{fig3}
 \end{center}
\end{figure}

\begin{figure}[ht!]
  \begin{center}
   \epsfig{file=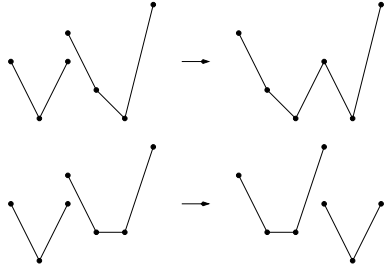,height=3.5cm}
   \caption{Transforming $a_{31}a_{13}a_{42}a_{21}a_{15}$ and $a_{31}a_{13}a_{42}a_{22}a_{25}$.}
   \label{fig4}
  \end{center}
\end{figure}

\begin{lemma} \label{cf2}
 If $A$ is a Cartier-Foata matrix, $C$ is a right-quantum matrix.
\end{lemma}
\begin{prf}
 Choose $i,j,k > n$, $i \neq j$. The product
 $c_{ik}c_{jk}$ is the sum of terms of the form
 $$a_{ii_1}a_{i_1i_2}\cdots a_{i_p k} a_{jj_1}a_{j_1j_2}\cdots a_{j_r k}$$
 for $p,r \geq 0$, $i_1,\ldots,i_p,j_1,\ldots,j_r \leq n$. Note that with the (possible) exception of $i,j,k$, all other terms appear as starting heights exactly as many times as they appear as ending heights.\\
 Identify this term with a sequence of steps, as described in Section \ref{alg}. We will perform a series of switches of steps that will transform such a term into a term of $c_{jk}c_{ik}$.\\
 The variable $a_{jj_1}$ (or $a_{jk}$ if $r = 0$) commutes with all variables that appear before it. In other words, in the algebra $\p A$, the expressions
 $$a_{ii_1}a_{i_1i_2}\cdots a_{i_p k} a_{jj_1}a_{j_1j_2}\cdots a_{j_r k}$$
 and
 $$a_{jj_1}a_{ii_1}a_{i_1i_2}\cdots a_{i_p k} a_{j_1j_2}\cdots a_{j_r k}$$
 are the same modulo the ideal $\p I_{\cf}$ generated by $a_{ik}a_{jl} - a_{jl}a_{ik}$ for $i \neq j$. Graphically, we can keep switching the step $j \to j_1$ with the step to its left until it is at the beginning of the sequence.\\
 If $r=0$, we are already done. If not, take the first step to the right of $a_{jj_1}$ that has starting height $j_1$; such a step certainly exists -- for example $j_1 \to j_2$. Without changing the expression modulo $\p I_{\cf}$, we can switch this step with the ones to the left until it is just right of $j \to j_1$. Continue this procedure; eventually, our sequence will have been transformed into an expression 
 of the form
 $$a_{jj_1'}a_{j_1'j_2'}\cdots a_{j_{r'}'k}a_{ii_1'}a_{i_1'i_2'}\cdots a_{i_{p'}' k}$$
 which will be equal modulo $\p I_{\cf}$ to the expression we started with.\\
 As an example, take $m=5$, $n=2$, $i=3$, $j=5$, $k=4$ and the term $a_{31}a_{12}a_{24}a_{52}a_{22}a_{24}$. The steps shown in Figure \ref{fig3} transform it into $a_{52}a_{24}a_{31}a_{12}a_{22}a_{24}$.\\
 It is clear that applying the same procedure to the result, but with the roles of $i$'s and $j$'s interchanged, gives the original sequence. This proves that indeed $c_{ik}c_{jk}=c_{jk}c_{ik}$.\\
 The proof of the other relation \eqref{cf8} is similar and we will only sketch it. Choose $i,j,k,l>n$, $i \neq j$, $k \neq l$. Then $c_{ik}c_{jl} + c_{il}c_{jk}$ is the sum of terms of the form
 $$a_{ii_1}a_{i_1i_2}\cdots a_{i_p k} a_{jj_1}a_{j_1j_2}\cdots a_{j_r l}$$
 and of the form
 $$a_{ii_1}a_{i_1i_2}\cdots a_{i_p l} a_{jj_1}a_{j_1j_2}\cdots a_{j_r k}$$
 for $p,r \geq 0$, $i_1,\ldots,i_p,j_1,\ldots,j_r \leq n$. Applying the same procedure as above to the first term yields either
 $$a_{jj_1'}a_{j_1'j_2'}\cdots a_{j_{r'}'k}a_{ii_1'}a_{i_1'i_2'}\cdots a_{i_{p'}' l}$$
 or 
 $$a_{jj_1'}a_{j_1'j_2'}\cdots a_{j_{r'}'l}a_{ii_1'}a_{i_1'i_2'}\cdots a_{i_{p'}' k},$$
 this procedure is reversible and it yields the desired identity. 
 See Figure \ref{fig4} for examples with $m=5$, $n=2$, $i=3$, $j=4$, $k=3$, $l=5$.\qed
\end{prf}

If $A$ is Cartier-Foata, Proposition \ref{cf1} implies
$$(I-A)^{-1}_{n+1,n+1}(I-A^{n+1,n+1})^{-1}_{n+2,n+2}\cdots =\det{}^{-1}(I-A) \cdot \det(I-A_0).$$
By Lemma \ref{cf2}, $C$ is right-quantum, so by Proposition \ref{cf1}
$$(I-C)^{-1}_{n+1,n+1}(I-C^{n+1,n+1})^{-1}_{n+2,n+2}\cdots=\det{}^{-1}(I-C),$$
and hence
$$\det{}^{-1}(I-A_0) \cdot \det(I-A)=\det(I-C).$$
In the classical Sylvester's identity, the entries of $I-C$ are also expressed as determinants. The following is an analogue of Proposition \ref{comm2}.

\begin{prop} \label{cf5}
 If $A$ is Cartier-Foata, then
 \begin{equation} \label{cf4}
  c_{ij} = - \det{}^{-1}(I-A_0) \cdot \det \begin{pmatrix} I-A_0 & -a_{*j} \\ -a_{i*} & -a_{ij} \end{pmatrix}.
 \end{equation}
\end{prop}
\begin{prf}
 We can repeat the proof of Proposition \ref{comm2} almost verbatim. We have
 $$(1-c_{ij})^{-1} = \left(\left( I - \begin{pmatrix} A_0 & a_{*j} \\ a_{i*} & a_{ij} \end{pmatrix}\right)^{-1}\right)_{ij},$$
 and because the matrix
 $$\begin{pmatrix} A_0 & a_{*j} \\ a_{i*} & a_{ij} \end{pmatrix}$$
 is still Cartier-Foata, Proposition \ref{cf1} shows that this is equal to
 $$\det{}^{-1} \left( I - \begin{pmatrix} A_0 & a_{*j} \\ a_{i*} & a_{ij} \end{pmatrix} \right) \cdot \det(I-A_0).$$
 We get
 $$1-c_{ij} = \det{}^{-1}(I-A_0)\cdot \det \left( I - \begin{pmatrix} A_0 & a_{*j} \\ a_{i*} & a_{ij} \end{pmatrix} \right)  = $$
 $$ = \det{}^{-1}(I-A_0)\cdot \left( \det \begin{pmatrix} I - A_0 & -a_{*j} \\ -a_{i*} & -a_{ij} \end{pmatrix} + \det \begin{pmatrix} I - A_0 & 0 \\ -a_{i*} & 1 \end{pmatrix} \right) =$$
 $$ = \det{}^{-1}(I-A_0)\cdot \left( \det \begin{pmatrix} I - A_0 & -a_{*j} \\ -a_{i*} & -a_{ij} \end{pmatrix} + \det (I - A_0) \right)= $$
 $$= \det{}^{-1}(I-A_0)\cdot \det \begin{pmatrix} I - A_0 & -a_{*j} \\ -a_{i*} & -a_{ij} \end{pmatrix} + 1. \eqno \qed$$
\end{prf}

We have proved the following.

\begin{thm}[Cartier-Foata Sylvester's identity] \label{cf3}
 Let $A=(a_{ij})_{m \times m}$ be a Cartier-Foata matrix, and choose $n < m$. Let $A_0,a_{i*},a_{*j}$ be defined as above, and let
 $$c_{ij} = - \det{}^{-1}(I-A_0) \cdot \det \begin{pmatrix} I-A_0 & -a_{*j} \\ -a_{i*} & -a_{ij} \end{pmatrix}, \quad C=(c_{ij})_{n+1 \leq i,j \leq m}.$$
 Then
 $$\det{}^{-1}(I-A_0) \cdot \det(I-A)=\det (I-C).\eqno \qed$$
\end{thm}

\section{Right-quantum analogue} \label{rq}

The right-quantum version of the Sylvester's identity is very similar; we will prove a right-quantum version of Lemma \ref{cf2} and Proposition \ref{cf5}, and a right-quantum version of Theorem \ref{cf3} will follow.

\medskip

The only challanging part is the following.

\begin{lemma} \label{rq1}
 If $A$ is a right-quantum matrix, so is $C$.
\end{lemma}
\begin{prf}
 Choose $i,j,k > n$, $i \neq j$. Instead of dealing directly with the equality $c_{ik}c_{jk}=c_{jk}c_{ik}$, we will prove an equivalent identity.\\
 Denote by $\p P_{ij}^k(k_1,k_2,\ldots,k_n)$ the set of sequences of $k_1+\ldots+k_n+2$ steps with the following properties:
 \begin{itemize}
  \item starting heights form a non-decreasing sequence;
  \item each $r$ between $1$ and $n$ appears exactly $k_r$ times as a starting height and exactly $k_r$ times as an ending height;
  \item $i$ and $j$ appear exactly once as starting heights;
  \item $k$ appears exactly twice as an ending height.
 \end{itemize}
 For $m=5,n=2,i=3,j=5,k=4,k_1=1,k_2=1$, all such sequences are shown in Figure \ref{fig5}.
 
 \begin{figure}[ht!]
  \begin{center}
   \includegraphics{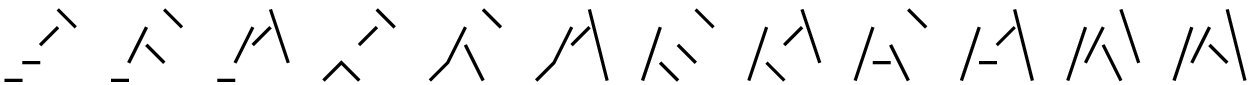}
   \caption{Sequences in the set $\p P_{35}^4(1,1)$.}
   \label{fig5}
  \end{center}
 \end{figure}
 
 We will do something very similar to the proof of Lemma \ref{cf2}: we will perform switches on sequences in $\p P_{ij}^k(k_1,k_2,\ldots,k_n)$ until they are transformed into sequences of the form $P_1P_2P_3$, where:
 \begin{itemize}
  \item $P_1$ is a path from $i$ to $k$ with all intermediate heights $\leq n$;
  \item $P_2$ is a path from $j$ to $k$ with all intermediate heights $\leq n$;
  \item $P_3$ is a sequence of steps with non-decreasing heights, with all heights $\leq n$, and with the number of steps with starting height $r$ equal to the number of steps with ending height $r$ for all $r$.
 \end{itemize}
 
 Namely, we move the step $i \to i'$ to the first place, the first step of the form $i' \to i''$ to the second place, etc. For example, the sequence $a_{11}a_{24}a_{34}a_{52}$ will be tranformed into $a_{34}a_{52}a_{24}a_{11}$, see Figure \ref{fig8}.
 
 \begin{figure}[ht]
  \begin{center}
   \includegraphics{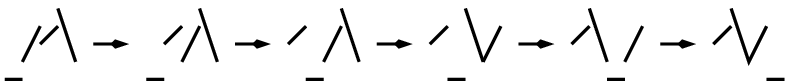}
  \end{center}
  \caption{Transforming $a_{11}a_{24}a_{34}a_{52}$ into $a_{34}a_{52}a_{24}a_{11}$.}
  \label{fig8}
 \end{figure}

 Of course, we have to prove that this can be done without changing the sum modulo the ideal $\p I_{\rtq}$ generated by relations \eqref{cf7}--\eqref{cf8}, and this is done in exactly the same way as the proof in \cite[Section 4]{konvalinka}. Figure \ref{fig6} is an example for $m=5,n=2,i=3,j=5,k=4,k_1=1,k_2=1$; each column corresponds to a transformation of an element of $\p P_{35}^4(1,1)$, if two elements in the same row have the same label, their sum can be transformed into the sum of the corresponding elements in the next row by use of the relation \eqref{cf8}, and if an element is not labeled it either means that it is transformed into the corresponding element in the next row by use of the relation \eqref{cf7} or is already in the required form. 
 
 \begin{figure}[ht]
  \begin{center}
   \epsfig{file=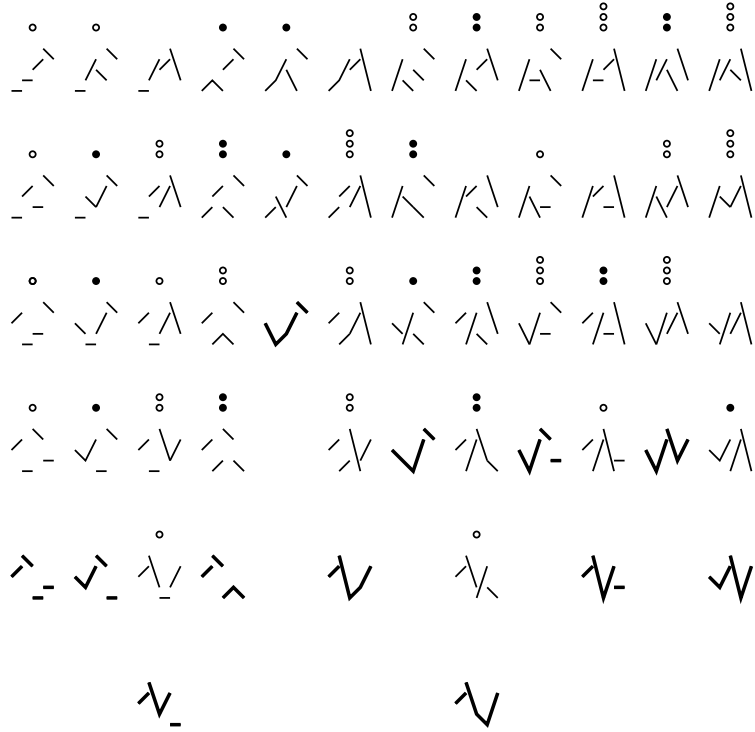,width=10cm}
  \end{center}
  \caption{Transforming the sequences in $\p P_{35}^4(1,1)$ into terms of $c_{34}c_{54}S$.}
  \label{fig6}
 \end{figure}
 
 This means that the sum of all elements of $\p P_{ij}^k(k_1,k_2,\ldots,k_n)$ over all $k_1,\ldots,k_n$ is modulo $\p I_{\rtq}$ equal to
 $$c_{ik}c_{jk} S,$$
 where $S$ is the sum over all sequences of steps with the following properties:
 
 \begin{itemize}
  \item starting heights form a non-decreasing sequence;
  \item starting and ending heights are all between $1$ and $n$;
  \item each $r$ between $1$ and $n$ appears as many times as a starting height as an ending height.
 \end{itemize}
  
 Of course, we can also reverse the roles of $i$ and $j$, and this proves that the sum of all elements of $\p P_{ij}^k(k_1,k_2,\ldots,k_n)$ is modulo $\p I_{\rtq}$ also equal to
 $$c_{jk}c_{ik}S,$$
 see Figure \ref{fig7}.
 
 \begin{figure}[ht]
  \begin{center}
   \epsfig{file=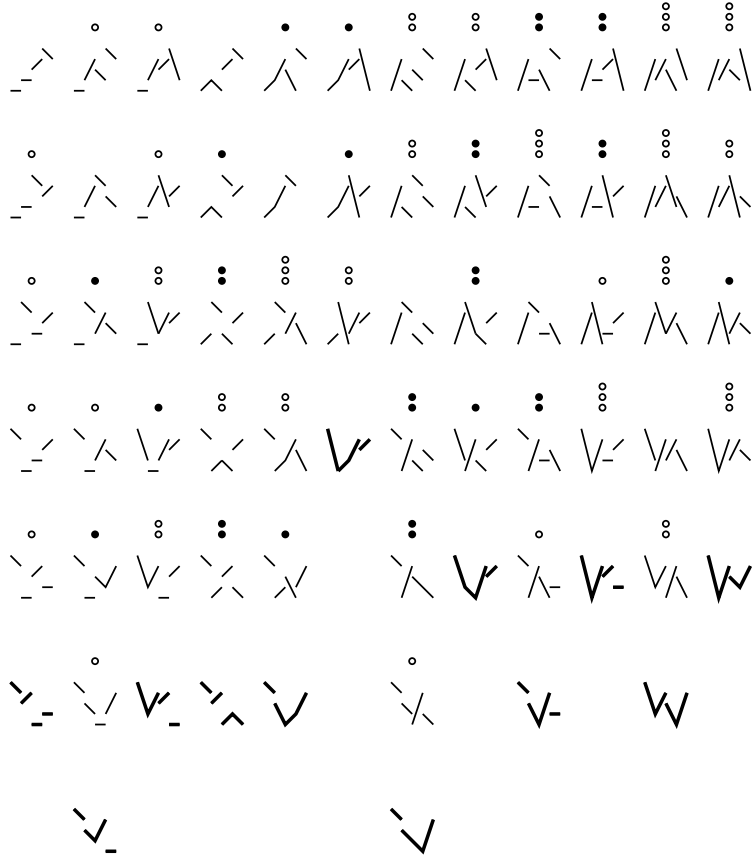,width=10cm}
  \end{center}
  \caption{Transforming the sequences in $\p P_{35}^4(1,1)$ into terms of $c_{54}c_{34}S$.}
  \label{fig7}
 \end{figure}
  
 Hence, modulo $\p I_{\rtq}$,
 \begin{equation} \label{rq5}
  c_{ik}c_{jk} S = c_{jk}c_{ik}S.
 \end{equation}
 But $S = 1 + a_{11}+\ldots+ a_{nn} + a_{11}a_{22}+a_{12}a_{21}+\ldots$ is an \emph{invertible} element of $\p A$, so \eqref{rq5} implies
 $$c_{ik}c_{jk} = c_{jk}c_{ik},$$
 provided $A$ is a right-quantum matrix.\\
 The proof of the other relation is almost completely analogous. Now we take $i \neq j$, $k \neq l$, and define $\p P_{ij}^{kl}(k_1,k_2,\ldots,k_n)$ as the set of sequences of $k_1+\ldots+k_n+2$ steps with the following properties:
 \begin{itemize}
  \item starting heights form a non-decreasing sequence;
  \item each $r$ between $1$ and $n$ appears exactly $k_r$ times as a starting height and exactly $k_r$ times as an ending height;
  \item $i$ and $j$ appear exactly once as starting heights;
  \item $k$ and $l$ appear exactly once as ending heights.
 \end{itemize}
 A similar reasoning shows that the sum over all elements of $\p P_{ij}^{kl}(k_1,k_2,\ldots,k_n)$ is equal both to $(c_{ik}c_{jl} + c_{il}c_{jk})S$ and to $(c_{jl}c_{ik} + c_{jk}c_{il})S$ modulo $\p I_{\rtq}$, which implies $c_{ik}c_{jl} + c_{il}c_{jk} = c_{jl}c_{ik} + c_{jk}c_{il}$.\qed
\end{prf}

\begin{prop} \label{rq2}
 If $A$ is right-quantum, then
 \begin{equation} \label{rq3}
  c_{ij} = -\det{}^{-1}(I-A_0) \cdot \det \begin{pmatrix} I-A_0 & -a_{*j} \\ -a_{i*} & -a_{ij} \end{pmatrix}.
 \end{equation}
\end{prop}
\begin{prf}
 The proof is exactly the same as the proof of Proposition \ref{cf5}.\qed
\end{prf}

\begin{thm}[right-quantum Sylvester's identity] \label{rq4}
 Let $A=(a_{ij})_{m \times m}$ be a right-quantum matrix, and choose $n < m$. Let $A_0,a_{i*},a_{*j}$ be defined as above, and let
  $$c_{ij} = - \det{}^{-1}(I-A_0) \cdot \det \begin{pmatrix} I-A_0 & -a_{*j} \\ -a_{i*} & -a_{ij} \end{pmatrix}, \quad C=(c_{ij})_{n+1 \leq i,j \leq m}.$$
  Then
 $$\det{}^{-1}(I-A_0) \cdot \det(I-A)=\det (I-C).\eqno \qed$$
\end{thm}

\section{$q$-Cartier-Foata analogue} \label{cfq}

Let us find a quantum extension of Theorem \ref{cf3}. Fix $q \in \C \setminus \{0\}$. We say that a matrix $A = (a_{ij})_{m \times m}$ is $q$-Cartier-Foata if 
\begin{eqnarray}
 a_{jl}  a_{ik}   &  = &   a_{ik}  a_{jl}  \ \
 \text{for}  \ \ i < j,   k < l , \label{cfq4}
  \\
 a_{jl}  a_{ik}   & = &   q^2   a_{ik}   a_{jl} 
 \ \ \text{for}  \ \ i < j,   k > l ,  \label{cfq5} \\
 a_{jk}  a_{ik}   & = &   q   a_{ik}  a_{jk} 
 \ \ \text{for}  \ \ i < j \label{cfq6},
\end{eqnarray}
and $q$-right-quantum if
\begin{eqnarray}
 a_{jk}  a_{ik} & = & q   a_{ik}a_{jk}  \ \
 \text{for all}  \ \ i < j ,
 \label{cfq7} \\
 a_{ik}  a_{jl}  -  q^{-1} a_{jk}  a_{il}   & = &  
 a_{jl}  a_{ik}  -q   a_{il}  a_{jk} 
    \ \ \text{for all}  \ \ i < j,k<l . \label{cfq8}
\end{eqnarray}

In the following two sections, the weight $w(\lambda,\mu)$ will be equal to $q^{\inv \mu - \inv \lambda}$. For example,
$$\detq (I-A) = \sum_{J \subseteq [m]} (-1)^{|J|} \detq A_J,$$
where
$$\detq A_J = \detq (a_{ij})_{i,j \in J} = \sum_{\sigma \in S_J} (-q)^{-\inv \sigma}a_{\sigma(j_1)j_1} \cdots a_{\sigma(j_k)j_k}$$
for $J=\set{j_1<j_2<\ldots<j_k}$.

\medskip

The following extends Proposition \ref{cf1}. See Section \ref{final} for a discussion of the proof.

\begin{prop} \label{cfq1}
 If $A=(a_{ij})_{m \times m}$ is a $q$-Cartier-Foata or a $q$-right-quantum matrix, we have
 $$\left(\frac1{I-A_{[ij]}}\right)_{ij} \, = \, (-1)^{i+j} \frac1{\detq(I-A)} \, \cdot \,
   \detq\left(I-A\right)^{ji}$$
 for all $i,j$, where
 $$A_{[ij]} = \begin{pmatrix} q^{-1} a_{11} & \cdots & q^{-1} a_{1j} & a_{1,j+1} & \cdots & a_{1m} \\ \vdots & \ddots & \vdots & \vdots & \ddots & \vdots \\ q^{-1} a_{i-1,1} & \cdots & q^{-1} a_{i-1,j} & a_{i-1,j+1} & \cdots & a_{i-1,m} \\ a_{i1} & \cdots & a_{ij} & q a_{i,j+1} & \cdots & q a_{i,m} \\ \vdots & \ddots & \vdots & \vdots & \ddots & \vdots \\ a_{m1} & \cdots & a_{mj} & q a_{m,j+1} & \cdots & q a_{mm} \end{pmatrix}. \eqno \qed$$
\end{prop}

We will use Theorem \ref{sylv2} for the matrix $A_{[ij]}$. Let us find the corresponding $C=(c'_{i'j'})_{n+1 \leq i',j' \leq m}$. Denote
$$a_{i'j'} + q^{-1} a_{i'*} (I- q^{-1} A_0)^{-1} a_{*j'}$$
by $c_{i'j'}$ for $i',j' > n$. If $i' < i,j' \leq j$, we have
$$c'_{i'j'}=q^{-1} a_{i'j'} + (q^{-1} a_{i'*}) (I- q^{-1} A_0)^{-1} (q^{-1} a_{*j'}) = q^{-1} c_{i'j'};$$
if $i' < i, j' > j$, we have
$$c'_{i'j'}=a_{i'j'} + (q^{-1} a_{i'*}) (I- q^{-1} A_0)^{-1} a_{*j'} = c_{i'j'};$$
if $i' \geq i, j' \leq j$, we have
$$c'_{i'j'}=a_{i'j'} + a_{i'*} (I- q^{-1} A_0)^{-1} (q^{-1} a_{*j'}) = c_{i'j'};$$
and if $i' \geq i, j' > j$, we have
$$c'_{i'j'}=q a_{i'j'} + a_{i'*} (I- q^{-1} A_0)^{-1} a_{*j'} = q c_{i'j'}.$$

We have proved the following.

\begin{prop} \label{cf9}
 With $A_{[ij]}$ as above and with $C=(c_{i'j'})_{n+1 \leq i',j' \leq m}$ for
 $$c_{i'j'} = a_{i'j'} + a_{i'*} (I-q^{-1}A_0)^{-1} (q^{-1}a_{*j'}),$$
 we have
 $$(I-A_{[ij]})^{-1}_{i'j'} = (I-C_{[ij]})^{-1}_{i'j'}. \eqno \qed$$
\end{prop}

\begin{rmk} \label{cf10}
 Let us present a slightly different proof of the proposition. Another way to characterize $A_{[ij]}$ is to say that the entry $a_{kl}$ has weight $q$ to the power of
 $${\left\{ \begin{matrix} 1 \colon l > j \\ 0 \colon l \leq j \end{matrix} \right. \: - \: \left\{ \begin{matrix} 1 \colon k < i \\ 0 \colon k \geq i \end{matrix} \right.} \: .$$
 That means that in $\left(A_{[ij]}^\ell\right)_{i_1i_\ell}$,
 $$a_{i_1i_2}a_{i_2i_3}\cdots a_{i_{\ell-1} i_\ell}$$
 has weight
 $$q^{|\set{r \colon i_r> j}|-|\set{r \colon i_r < i}|}.$$
 Assume that we have a decomposition of a path of length $\ell$ from $i'$ to $j'$, $i',j' > n$, as in Section \ref{sylv}, say $a_{\lambda,\mu} = a_{i'\lambda_1,\lambda_1i_1}a_{i_1\lambda_2,\lambda_2i_2}\cdots a_{i_{p-1}\lambda_p,\lambda_p j'}$, with all elements of $\lambda_r$ at most $n$, $i_r > n$, and the length of $\lambda_r$ equal to $\ell_r$. Put $i_0 = i'$, $i_{p+1}=j'$.The number of indices of $\lambda = i'\lambda_1\ldots \lambda_p$ that are strictly smaller than $i$ is clearly 
 $$ \sum_{r=1}^p \ell_r + |\set{r \colon i_r < i}| = \ell - p + |\set{r \colon i_r < i}|,$$
 and the number of indices of $\mu = \lambda_1\ldots \lambda_p j'$ that are strictly greater than $j$ is  $|\set{r \colon i_r > j}|$. Therefore the path $a_{\lambda,\mu}$ is weighted by
 $$ q^{- \ell + p + |\set{r \colon i_r > j}| - |\set{r \colon i_r < i}|}.$$
 On the other hand, take a term $a_{\lambda,\mu}=a_{i'\lambda_1,\lambda_1i_1}a_{i_1\lambda_2,\lambda_2i_2}\cdots a_{i_{p-1}\lambda_p,\lambda_p j'}$ (with $\lambda_r,i_r,\ell_r$ as before) of $(C_{[ij]}^\ell)_{i'j'}$. Each $a_{i_{r-1}\lambda_r,\lambda_ri_r}$ has weight $q^{-\ell_r}$ as an element of $C$, and $a_{\lambda,\mu}$ has the additional weight
 $$q^{|\set{r \colon i_r > j}| - |\set{r \colon i_r < i}|}$$
 as a term of $(C_{[ij]}^\ell)_{i'j'}$. The proposition follows. \qed
\end{rmk}

In what follows, the crucial observation will be the following. Take $a_{\lambda,\mu}$, $\lambda = \lambda_1 ij \lambda_2$, $\mu = \mu_1 kl \mu_2$, $\lambda'=\lambda_1 ji \lambda_2$, $\mu' = \mu_1 lk \mu_2$ for $i < j$. Then
$$q^{\inv \mu - \inv \lambda} a_{\lambda,\mu} = q^{\inv \mu' - \inv \lambda'} a_{\lambda' \mu'} \  \mod \  \p I_{q-\cf},$$
where $\p I_{q-\cf}$ is the ideal of $\p A$ generated by the equations \eqref{cfq4}--\eqref{cfq6}.

\medskip

We show this by considering in turn each of the following possibilities:
\begin{enumerate}
 \item \label{cfq3} $i < j, k < l$
 \item $i < j, k > l$
 \item $i < j, k = l$
\end{enumerate}

For example, to prove case \eqref{cfq3}, note that $a_{jl} a_{ik} - a_{ik} a_{jl}$ is a generator of $\p I_{q-\cf}$, and that $\inv \mu' = \inv \mu + 1$ and $\inv \lambda' = \inv \lambda + 1$. Other cases are similarly straightforward.	

\begin{lemma} \label{cfq2}
 If $A$ is a $q$-Cartier-Foata matrix, $C$ is a $q$-right-quantum matrix.
\end{lemma}
\begin{prf}
 We will adapt the proof of Lemma \ref{cf2}. Choose $i,j,k > n$, $i < j$. The product
 $c_{ik}c_{jk}$ is the sum of terms of the form
 $$q^{-p-r} a_{ii_1}a_{i_1i_2}\cdots a_{i_p k} a_{jj_1}a_{j_1j_2}\cdots a_{j_r k}$$
 for $p,r \geq 0$, $i_1,\ldots,i_p,j_1,\ldots,j_r \leq n$.\\
 Without changing the expression modulo $\p I_{q-\cf}$, we can repeat the procedure in the proof of Lemma \ref{cf2}, keeping track of weight changes. The resulting expression
 $$a_{jj_1'}a_{j_1'j_2'}\cdots a_{j_{r'}'k}a_{ii_1'}a_{i_1'i_2'}\cdots a_{i_{p'}' k}$$
 will, by the discussion preceding the lemma, have weight $q^{-1-r'-p'}$ (the extra $-1$ comes from the fact that the step with starting height $j$ is now to the left of the step with starting height $i$), In other words,
 $$c_{jk}c_{ik} = q c_{ik} c_{jk}.$$
 The proof of the other relation is completely analogous.\qed
\end{prf}

If $A$ is $q$-Cartier-Foata, Proposition \ref{cfq1} implies
$$(I-A_{[n+1,n+1]})^{-1}_{n+1,n+1}(I-\left(A^{n+1,n+1}\right)_{[n+2,n+2]})^{-1}_{n+2,n+2}\cdots =\detq{}^{-1}(I-A) \cdot \detq(I-A_0).$$
By Lemma \ref{cfq2}, $C$ is $q$-right-quantum, so by Proposition \ref{cfq1}
$$(I-C_{[n+1,n+1]})^{-1}_{n+1,n+1}(I-\left(C^{n+1,n+1}\right)_{[n+2,n+2]})^{-1}_{n+2,n+2}\cdots=\detq{}^{-1}(I-C),$$
and hence
$$\detq{}^{-1}(I-A_0) \cdot \detq(I-A)=\detq(I-C).$$

\medskip

The final step is to write entries of $C$ as quotients of quantum determinants.

\begin{prop} \label{cfq10}
 If $A$ is $q$-Cartier-Foata, then
 $$c_{ij} = -\detq{}^{-1}(I-A_0) \cdot \detq \begin{pmatrix} I-A_0 & -a_{*j} \\ -a_{i*} & -a_{ij} \end{pmatrix}.$$
\end{prop}
\begin{prf}
 Again,
 $$(1-c_{ij})^{-1} = \left(\left( I - \begin{pmatrix} q^{-1} A_0 & q^{-1} a_{*j} \\ a_{i*} & a_{ij} \end{pmatrix}\right)^{-1}\right)_{ij},$$
 and because the matrix
 $$\begin{pmatrix} A_0 & a_{*j} \\ a_{i*} & a_{ij} \end{pmatrix}$$
 is still $q$-Cartier-Foata, Proposition \ref{cfq1} shows that this is equal to
 $$\detq{}^{-1} \left( I - \begin{pmatrix} A_0 & a_{*j} \\ a_{i*} & a_{ij} \end{pmatrix} \right) \cdot \detq(I-A_0).$$
 The rest of the proof is exactly the same as in Proposition \ref{cf5}, with $\detq$ playing the role of $\det$.\qed
\end{prf}

We have proved the following.

\begin{thm}[$q$-Cartier-Foata Sylvester's identity] \label{cfq9}
 Let $A=(a_{ij})_{m \times m}$ be a $q$-Cartier-Foata matrix, and choose $n < m$. Let $A_0,a_{i*},a_{*j}$ be defined as above, and let
 $$c_{ij}^q = -\detq\!{}^{-1}(I-A_0) \cdot \detq \begin{pmatrix} I-A_0 & -a_{*j} \\ -a_{i*} & -a_{ij} \end{pmatrix}, \quad C^q = (c_{ij}^q)_{n+1 \leq i,j \leq m}.$$
 Then
 $$\detq\!{}^{-1}(I-A_0) \cdot \detq(I-A)=\detq (I-C^q). \eqno \qed$$
\end{thm}
 
\section{$q$-right-quantum analogue} \label{qrq}

The results of the previous two sections easily extend to a $q$-right-quantum Sylvester's identity. Denote the ideal generated by relations \eqref{cfq7}--\eqref{cfq8} by $\p I_{q-\rtq}$. It is easy to see that if $\lambda = \lambda_1 ij \lambda_2$, $\mu = \mu_1 kl \mu_2$, $\lambda'=\lambda_1 ji \lambda_2$, $\mu' = \mu_1 lk \mu_2$ and if $i < j$, then
$$q^{\inv \mu - \inv \lambda} a_{\lambda,\mu} + q^{\inv \mu' - \inv \lambda} a_{\lambda,\mu'}= q^{\inv \mu - \inv \lambda'} a_{\lambda', \mu} + q^{\inv \mu' - \inv \lambda'} a_{\lambda',\mu'} \  \mod \  \p I_{q-\rtq}.$$

\begin{lemma} \label{qrq1}
 If $A$ is a $q$-right-quantum matrix, so is $C$.
\end{lemma}
\begin{prf}
 This is a weighted analogue of Lemma \ref{rq1}. The sum over all elements of $\p P_{ij}^k(k_1,\ldots,k_n)$ with $a_{\lambda, \mu}$ weighted by $q^{\inv \mu - \inv \lambda}=q^{\inv \mu}$ is modulo $\p I_{q-\rtq}$ equal to both $c_{ik}c_{jk}S$ and $q^{-1} c_{jk}c_{ik}S$; this implies the relation \eqref{cfq7} for elements of $C$, and the proof of \eqref{cfq8} is completely analogous. \qed
\end{prf}

\begin{prop} \label{qrq2}
 If $A$ is $q$-right-quantum, then
 $$c_{ij} = -\detq{}^{-1}(I-A_0) \cdot \detq \begin{pmatrix} I-A_0 & -a_{*j} \\ -a_{i*} & -a_{ij} \end{pmatrix}.$$
\end{prop}
\begin{prf}
 The proof is exactly the same as the proof of Proposition \ref{cfq10}.\qed
\end{prf}

Proposition \ref{cf9}, Lemma \ref{qrq1} and Proposition \ref{qrq2} imply the following theorem.

\begin{thm}[$q$-right-quantum Sylvester's identity]
 Let $A=(a_{ij})_{m \times m}$ be a $q$-right-quantum matrix, and choose $n < m$. Let $A_0,a_{i*},a_{*j}$ be defined as above, and let
 $$c_{ij}^q = -\detq\!{}^{-1}(I-A_0) \cdot \detq \begin{pmatrix} I-A_0 & -a_{*j} \\ -a_{i*} & -a_{ij} \end{pmatrix}, \quad C^q = (c_{ij}^q)_{n+1 \leq i,j \leq m}.$$
 Then
 $$\detq\!{}^{-1}(I-A_0) \cdot \detq(I-A)=\detq (I-C^q). \eqno \qed$$
\end{thm}

\section{$q_{ij}$-Cartier-Foata analogue} \label{cfqij}

Now let us prove a multiparameter extension of Theorem \ref{cfq9}. Choose $q_{ij} \neq 0$ for $i<j$, and recall that a matrix $A = (a_{ij})_{m \times m}$ is $\s q$-Cartier-Foata if 
\begin{eqnarray}
 q_{kl}   a_{jl}  a_{ik}   & = &   q_{ij}   a_{ik}  a_{jl} 
 \ \ \text{for} \ \ i < j,k < l , \label{cfqij2} \\
 a_{jl}  a_{ik}   & = &   q_{ij}   q_{lk}   a_{ik}  a_{jl} 
 \ \ \text{for}
 \ \ i < j,k > l ,
 \label{cfqij3} \\
 a_{jk}  a_{ik}  & = &  q_{ij}   a_{ik}  a_{jk} 
 \ \ \text{for}  \ \ i < j ,
 \label{cfqij4}
\end{eqnarray}
and $\s q$-right-quantum if
\begin{eqnarray}
a_{jk}  a_{ik} &  = &  q_{ij}   a_{ik}  a_{jk}  \ \
\text{for all} \ \ i < j, \label{cfqij5} \\
a_{ik}  a_{jl}   -   q_{ij}^{-1}   a_{jk}  a_{il} &  = & 
q_{kl}q_{ij}^{-1}   a_{jl}  a_{ik}   -  
q_{kl}   a_{il}  a_{jk} \ \
\text{for all} \ \ i < j, \ k <l. \label{cfqij6}
\end{eqnarray}

\medskip

If we define $q_{ii} = 1$ and $q_{ji} = q_{ij}^{-1}$ for $i < j$, we can write the conditions \eqref{cfqij2}--\eqref{cfqij4} more concisely as
\begin{equation} \label{cfqij7}
 q_{kl} a_{jl}  a_{ik}  =  q_{ij}  a_{ik}  a_{jl},
\end{equation}
for all $i,j,k,l$, $i \neq j$, and \eqref{cfqij5}--\eqref{cfqij6} as
\begin{equation} \label{cfqij8}
a_{ik} a_{jl}  -  q_{ij}^{-1}  a_{jk}  a_{il}  = 
q_{kl}q_{ij}^{-1}  a_{jl}  a_{ik}  - 
q_{kl}  a_{il}  a_{jk}
\end{equation}
for all $i,j,k,l$, $i \neq j$.

\medskip

In the following two sections, the weight $w(\lambda,\mu)$ will be equal to
$$\prod_{(i,j) \in \p I(\mu)} q_{\mu_j \mu_i} \prod_{(i,j) \in \p I(\lambda)} q_{\lambda_j \lambda_i}^{-1}.$$
For example,
$$\detq (I-A) = \sum_{J \subseteq [m]} (-1)^{|J|} \detqq A_J,$$
where
$$\detq A_J = \detq (a_{ij})_{i,j \in J} = \sum_{\sigma \in S_J} \left( \prod_{p<q: \sigma(p)>\sigma(q)} q_{\sigma(q)\sigma(p)}^{-1} \right) a_{\sigma(j_1)j_1} \cdots a_{\sigma(j_k)j_k}$$
for $J=\set{j_1<j_2<\ldots<j_k}$.

\medskip

The following extends Proposition \ref{cfq1}, see Section \ref{final}.

\begin{prop} \label{cfqij9}
 If $A=(a_{ij})_{m \times m}$ is a $\s q$-Cartier-Foata matrix or a $\s q$-right-quantum matrix, we have
 $$\left(\frac1{I-A_{[ij]}}\right)_{ij} \, = \, (-1)^{i+j} \frac1{\detqq(I-A)} \, \cdot \,
   \detqq\left(I-A\right)^{ji}$$
 for all $i,j$, where
 $$A_{[ij]} \!  = \! \begin{pmatrix} q_{1i}^{-1} a_{11} & \cdots & q_{1i}^{-1} a_{1j} & q_{1i}^{-1}q_{j,j+1}a_{1,j+1} & \cdots & q_{1i}^{-1}q_{jm}a_{1m} \\ \vdots & \ddots & \vdots & \vdots & \ddots & \vdots \\ q_{i-1,i}^{-1} a_{i-1,1} & \cdots & q_{i-1,i}^{-1} a_{i-1,j} & q_{i-1,i}^{-1}q_{j,j+1} a_{i-1,j+1} & \cdots & q_{i-1,i}^{-1} q_{jm} a_{i-1,m} \\ a_{i1} & \cdots & a_{ij} & q_{j,j+1} a_{i,j+1} & \cdots & q_{jm} a_{i,m} \\ \vdots & \ddots & \vdots & \vdots & \ddots & \vdots \\ a_{m1} & \cdots & a_{mj} & q_{j,j+1} a_{m,j+1} & \cdots & q_{jm} a_{mm} \end{pmatrix}.\eqno \!\!\!\! \qed$$
\end{prop}

Assume that $q_{ij}=q_{i'j'}$ for $i,i' \leq n$, $j,j' > n$; denote this value by $q$. We will use Theorem \ref{sylv2} for the matrix $A_{[ij]}$ and the corresponding $C=(c'_{i'j'})_{n+1 \leq i',j' \leq m}$. Define
$$c_{i'j'} = a_{i'j'} + q^{-1} a_{i'*} (I- q^{-1} A_0)^{-1} a_{*j'}$$
for $i',j' > n$. If $i' < i,j' \leq j$, we have
$$c'_{i'j'}=q_{i'i}^{-1} a_{i'j'} + (q_{i'i}^{-1} a_{i'*}) (I- q^{-1} A_0)^{-1} (q^{-1} a_{*j'}) = q_{i'i}^{-1} c_{i'j'};$$
if $i' < i, j' > j$, we have
$$c'_{i'j'}=q_{i'i}^{-1}q_{jj'} a_{i'j'} + (q_{i'i}^{-1} a_{i'*}) (I- q^{-1} A_0)^{-1} (q^{-1}q_{jj'}a_{*j'}) = q_{i'i}^{-1}q_{jj'}c_{i'j'};$$
if $i' \geq i, j' \leq j$, we have
$$c'_{i'j'}=a_{i'j'} + a_{i'*} (I- q^{-1} A_0)^{-1} (q^{-1} a_{*j'}) = c_{i'j'};$$
and if $i' \geq i, j' > j$, we have
$$c'_{i'j'}=q_{jj'} a_{i'j'} + a_{i'*} (I- q^{-1} A_0)^{-1} (q^{-1} q_{jj'}a_{*j'}) = q_{jj'} c_{i'j'}.$$

We have proved the following.

\begin{prop}
 With $A_{[ij]}$ as defined above and with $C=(c_{i'j'})_{n+1 \leq i',j' \leq m}$ for
 $$c_{i'j'} = a_{i'j'} + a_{i'*} (I- q^{-1} A_0)^{-1} (q^{-1} a_{*j'}),$$
 we have
 $$(I-A_{[ij]})^{-1}_{i'j'} = (I-C_{[ij]})^{-1}_{i'j'}. \eqno \qed$$
\end{prop}

\begin{rmk}
 Another way to characterize $A_{[ij]}$ is to say that the entry $a_{kl}$ has weight
 $${\left\{ \begin{matrix} q_{jl} \colon l > j \\ 1 \colon l \leq j \end{matrix} \right. \: \cdot \:\left\{ \begin{matrix} q_{ki}^{-1} \colon k < i \\ 1 \colon k \geq i \end{matrix} \right.} \: .$$
 That means that in $\left(A_{[ij]}^\ell\right)_{i_1i_\ell}$,
 $$a_{i_1i_2}a_{i_2i_3}\cdots a_{i_{\ell-1} i_\ell}$$
 has weight
 $$\prod_{i_r > j} q_{j i_r} \cdot \prod_{i_r < i} q_{i_r i}^{-1}.$$
 An alternative way to prove the proposition is analogous to the proof of Proposition \ref{cf9} outlined in Remark \ref{cf10}.\qed
\end{rmk}

If $a_{\lambda,\mu}$, $\lambda = \lambda_1 ij \lambda_2$, $\mu = \mu_1 kl \mu_2$, $\lambda'=\lambda_1 ji \lambda_2$, $\mu' = \mu_1 lk \mu_2$ and if $i < j$, then
$$\left(\prod_{(i,j) \in \p I(\mu)} q_{\mu_j \mu_i}
\prod_{(i,j) \in \p I(\lambda)} q_{\lambda_j \lambda_i}^{-1}\right) a_{\lambda,\mu} = \left(\prod_{(i,j) \in I(\mu')} q_{\mu_j' \mu_i'}
\prod_{(i,j) \in I(\lambda')} q_{\lambda_j' \lambda_i'}^{-1}\right) a_{\lambda' \mu'} \  \mod   \p I_{q-\cf},$$
where $\p I_{\s q-\cf}$ is the ideal of $\p A$ generated by the equations \eqref{cfqij2}--\eqref{cfqij4}.

\medskip

As in the $q$-Cartier-Foata case, we show this by considering in turn each of the possibilities $k<l$, $k>l$, $k=l$.

\begin{lemma} \label{cfqij1}
 If $A$ is a $\s q$-Cartier-Foata matrix, $C$ is a $\s q$-right-quantum matrix.
\end{lemma}
\begin{prf}
 We will adapt the proof of Lemma \ref{cfq2}. Choose $i,j,k > n$, $i < j$. The product
 $c_{ik}c_{jk}$ is the sum of terms of the form
 $$q^{-p-r} a_{ii_1}a_{i_1i_2}\cdots a_{i_p k} a_{jj_1}a_{j_1j_2}\cdots a_{j_r k}$$
 for $p,r \geq 0$, $i_1,\ldots,i_p,j_1,\ldots,j_r \leq n$.\\
 Note that since
 $$q^{-p-r} = q_{j_1k}\cdots q_{j_r k} q_{i_1i}^{-1}\cdots q_{i_pi}^{-1} q_{j_1i}^{-1}\cdots q_{j_ri}^{-1} q_{j_1j}^{-1}\cdots q_{j_rj}^{-1},$$
 the weight of $a_{ii_1}a_{i_1i_2}\cdots a_{i_p k} a_{jj_1}a_{j_1j_2}\cdots a_{j_r k}$ is of the form 
 $$\prod_{(i,j) \in \p I(\mu)} q_{\mu_j \mu_i}
 \prod_{(i,j) \in \p I(\lambda)} q_{\lambda_j \lambda_i}^{-1}$$
 for $\lambda = ii_1\ldots i_p j j_1 \ldots j_r$ and $\mu = i_1\ldots i_p k j_1 \ldots j_r k$.
 Without changing the expression modulo $\p I_{\s q-\cf}$, we can repeat the procedure in the proof of Lemma \ref{cf2}, but changing the weight at each switch. The resulting expression
 $$a_{jj_1'}a_{j_1'j_2'}\cdots a_{j_{r'}'k}a_{ii_1'}a_{i_1'i_2'}\cdots a_{i_{p'}' k}$$
 will, by the discussion preceding the lemma, have weight
 $$q_{i_1'k} \cdots q_{i_{p'}' k} q_{j_1' j}^{-1} \cdots q_{j_{r'}' j}^{-1} q_{i_1' j}^{-1} \cdots q_{i_{p'}' j}^{-1} q_{i_1' i}^{-1} \cdots q_{i_{p'}' i}^{-1} q_{ij}^{-1}= q^{-r'-p'} q_{ij}^{-1}$$
 (the extra $q_{ij}^{-1}$ comes from the fact that the step with starting height $j$ is now to the left of the step with starting height $i$), In other words,
 $$c_{jk}c_{ik} = q_{ij} c_{ik} c_{jk}.$$
 The proof of the other relation is completely analogous.\qed
\end{prf}

If $A$ is $\s q$-Cartier-Foata, Proposition \ref{cfqij9} implies
$$(I-A_{[n+1,n+1]})^{-1}_{n+1,n+1}(I-\left(A^{n+1,n+1}\right)_{[n+2,n+2]})^{-1}_{n+2,n+2}\cdots =\detqq\!\!\!{}^{-1}(I-A) \cdot \detqq(I-A_0).$$
By Lemma \ref{cfq2}, $C$ is $\s q$-right-quantum, so by Proposition \ref{cfqij9}
$$(I-C_{[n+1,n+1]})^{-1}_{n+1,n+1}(I-\left(C^{n+1,n+1}\right)_{[n+2,n+2]})^{-1}_{n+2,n+2}\cdots=\detqq\!\!\!{}^{-1}(I-C),$$
and hence
$$\detqq\!\!\!{}^{-1}(I-A_0) \cdot \detqq(I-A)=\detqq(I-C).$$

\medskip

So far, the extension to the multiparameter case has been straightforward. However, we need something extra for the proof of the analogue of Proposition \ref{cfq10} since the matrix
$$\begin{pmatrix} A_0 & a_{*j} \\ a_{i*} & a_{ij} \end{pmatrix}$$
is in general not $\s q$-Cartier-Foata. It turns out that a special case of the first statement of Proposition \ref{cfqij9} holds under slightly weaker conditions; see Section \ref{final}.

\begin{prop} \label{cfqij10}
 Assume that for a matrix $A=(a_{ij})_{m \times m}$, $A=(a_{ij})_{m \times (m-1)}$ is a $\s q$-Cartier-Foata or a $\s q$-right-quantum matrix. Then
 $$\left(\frac1{I-A_{[mm]}}\right)_{mm} \, = \, \frac1{\detqq(I-A)} \, \cdot \,
   \detqq\left(I-A^{mm}\right),$$
 where $A_{[mm]}$ is defined as in Proposition \ref{cfqij9}.\qed
\end{prop}

In other words, even though only the first $n$ columns of 
$$\begin{pmatrix} A_0 & a_{*j} \\ a_{i*} & a_{ij} \end{pmatrix}$$
satisfy the $\s q$-Cartier-Foata condition, we still have
$$(1-c_{ij})^{-1} = \left(\left( I - \begin{pmatrix} q^{-1} A_0 & q^{-1} a_{*j} \\ a_{i*} & a_{ij} \end{pmatrix}\right)^{-1}\right)_{ij} = $$
$$ = \detqq\!\!\!{}^{-1} \left( I - \begin{pmatrix} A_0 & a_{*j} \\ a_{i*} & a_{ij} \end{pmatrix} \right) \cdot \detqq(I-A_0).$$

\begin{prop} \label{cfqij12}
 If $A$ is $\s q$-Cartier-Foata, then
 $$c_{ij} = -\detqq\!\!\!{}^{-1}(I-A_0) \cdot \detqq \begin{pmatrix} I-A_0 & -a_{*j} \\ -a_{i*} & -a_{ij} \end{pmatrix}. \eqno \qed$$
\end{prop}
\begin{prf}
 This follows from the previous proposition, using the same technique as in the proof of Proposition \ref{cf5}.\qed
\end{prf}

We have proved the following.

\begin{thm}[$\s q$-Cartier-Foata Sylvester's theorem]
 Let $A=(a_{ij})_{m \times m}$ be a $\s q$-Cartier-Foata matrix, and choose $n < m$. Let $A_0,a_{i*},a_{*j}$ be defined as above, and let
  $$c^{\s q}_{ij} = - \detqq\!\!\!{}^{-1}(I-A_0) \cdot \detqq \begin{pmatrix} I-A_0 & -a_{*j} \\ -a_{i*} & -a_{ij} \end{pmatrix}, \quad C^{\s q}=(c_{ij}^{\s q})_{n+1 \leq i,j \leq m}.$$
  Suppose $q_{ij}=q_{i'j'}$ for all $i,i' \leq n$ and $j,j' > n$. Then
 $$\detqq\!\!\!{}^{-1}(I-A_0) \cdot \detqq(I-A)=\detqq (I-C^{\s q}).$$
\end{thm}

\begin{rmk}
 It is important to note that the determinant $\detqq (I-C^{\s q})$ is with respect to $\p C$, the algebra generated by $c_{ij}$'s, not with respect to $\p A$. For example, for $n=2$ and $m=4$, we have
 $$\detqq (I-C^{\s q})=1-c_{33}^{\s q}-c_{44}^{\s q}+c_{33}^{\s q}c_{44}^{\s q}-q_{34}^{-1}c_{43}^{\s q}c_{34}^{\s q}.$$
\end{rmk}

\begin{rmk} \label{cfqij11}
 The condition $q_{ij}=q_{i'j'}$ whenever $i,i' \leq n$, $j,j' > n$ is indeed necessary, as shown by the following. Take $n=1$ and $m=3$. In $\detqq\!\!\!{}^{-1}(I-A_0) \cdot \detqq(I-A)$ we have the term
 $$-q_{12}^{-1}q_{13}^{-1}a_{21}a_{32}a_{13},$$
 while in $\detqq (I-C^{\s q})$ we have
 $$-q_{23}^{-1}(-a_{32})(-q_{12}^{-1}a_{21}a_{13})=-q_{12}^{-2}a_{21}a_{32}a_{13}. \eqno \qed$$ 
\end{rmk}
 
\section{$q_{ij}$-right-quantum analogue} \label{qrqij}

The results in this section are almost complete copies of proofs above.

\medskip

Assume we have a $\s q$-right-quantum matrix, with $q_{ij}=q$ for $i \leq n$, $j > n$. In the notation of the previous section, we have the following.

\begin{lemma}
 If $A$ is a $\s q$-Cartier-Foata matrix, $C$ is a $\s q$-right-quantum matrix.
\end{lemma}
\begin{prf}
 We use a combination of proofs of Lemmas \ref{qrq1} and \ref{cfqij1}.\qed
\end{prf}

\begin{prop}
 If $A$ is $\s q$-right-quantum, then
 $$c_{ij} = -\detqq\!\!\!{}^{-1}(I-A_0) \cdot \detqq \begin{pmatrix} I-A_0 & -a_{*j} \\ -a_{i*} & -a_{ij} \end{pmatrix}.$$
\end{prop}
\begin{prf}
 We use the same technique as in the proof of Proposition \ref{cfqij12}.\qed
\end{prf}

This finishes the proof of Theorem \ref{intro3}.

\section{The $\beta$-extension} \label{beta}

Theorem \ref{intro2} trivially implies that
$$(\det B)^\beta = (\det A)^\beta \cdot (\det A_0)^{\beta(m-n-1)}$$
for any $\beta \in \C$, where $a_{ij}$ are commutative variables and
$$b_{ij} = \det \begin{pmatrix} A_0 & a_{*j} \\ a_{i*} & a_{ij} \end{pmatrix}, \quad B = (b_{ij})_{n+1 \leq i,j \leq m}.$$

It is not immediately clear what the non-commutative extension of this could be. Of course, Theorem \ref{cf3} implies that
$$(\det (I-C))^\beta = \left(\det{}^{-1}(I-A_0) \cdot \det(I-A)\right)^\beta$$
for
$$c_{ij} = - \det{}^{-1}(I-A_0) \cdot \det \begin{pmatrix} I-A_0 & -a_{*j} \\ -a_{i*} & -a_{ij} \end{pmatrix}, \quad C=(c_{ij})_{n+1 \leq i,j \leq m},$$
where $A$ is a Cartier-Foata or right-quantum matrix, but this does not really tell us much about $(\det (I-C))^\beta$, especially when $\beta$ is not an integer. However, a technique similar to the proof of the $\beta$-extension of the non-commutative MacMahon Master Theorem, \cite[\S 10]{konvalinka}, gives a reasonable interpretation of $(\det (I-C))^\beta$ for $\beta \in \C$ when $A$ is a Cartier-Foata matrix.

\medskip

We will need some terminology from \cite{konvalinka}. A \emph{balanced sequence} (\emph{b-sequence}) is
a finite sequence of steps such that the number of steps starting at height~$i$ is equal to the number of steps ending at height~$i$, for all~$i$. We denote this number by $k_i$, and call $(k_1,\ldots,k_m)$ the \emph{type} of the b-sequence. An \emph{ordered sequence} (\emph{o-sequence}) is a b-sequence where the steps starting at smaller height always precede steps starting at larger heights. In other words, an o-sequence of type $(k_1,\ldots,k_m)$ is a sequence of $k_1$ steps starting at height~$1$, then $k_2$ steps starting at
height~$2$, etc., so that $k_i$ steps end at height $i$. Denote by $\s O(k_1,\ldots,k_m)$ the set of all o-sequences of type $(k_1,\ldots,k_m)$. Finally, consider a lattice path from $(0,1)$ to $(x_1,1)$ that never goes below $y=1$ or above $y=m$, then a lattice path from $(x_1,2)$ to $(x_2,2)$ that never goes below $y=2$ or above $y=m$, etc.; in the end, take a straight path from $(x_{m-1},m)$ to $(x_m,m)$. We will call this a \emph{path sequence} (\emph{p-sequence}). Observe that every p-sequence is also a b-sequence. Denote by
$\s P(k_1,\ldots,k_m)$ the set of all p-sequences of type $(k_1,\ldots,k_m)$.

\medskip

In \cite[\S 2]{konvalinka} a bijection 
$$\varphi \, : \, \s O(k_1,\ldots,k_m) \, \longrightarrow \, \s P(k_1,\ldots,k_m)$$
was defined (which proved various forms of the MacMahon Master Theorem) as follows. Take an o-sequence $\alpha$, and let $[0,x]$ be the maximal interval on which it is part of a p-sequence, i.e.\hspace{-0.07cm} the maximal interval $[0,x]$ on which the o-sequence has the property that if a step ends at level $i$, and the following step starts at level $j > i$, the o-sequence stays on or above height $j$ afterwards. Let~$i$ be the height at~$x$. Choose the step $(x',i) \to (x'+1,i')$ in the o-sequence that is the first to the right of $x$ that starts at level $i$ (such a step exists because an o-sequence is a balanced sequence). Continue switching this step with the one to the left until it becomes the step $(x,i) \to (x+1,i')$. The new object is part of a p-sequence at least on the interval $[0,x+1]$. Continuing this procedure we get a p-sequence $\vp(\alpha)$.

\medskip

A lattice path from $i$ to $i$ with each height appearing at most once as the starting height will be called a \emph{disjoint cycle}.

\medskip

For an o-sequence $a_{\lambda,\mu}$, take the corresponding p-sequence $a_{\lambda',\mu'}=\vp(a_{\lambda,\mu})$. If the first repeated height in $a_{\lambda',\mu'}$ is the starting height of the sequence, the sequence starts with a disjoint cycle; remove it and repeat the algorithm. If the first repeated height in $a_{\lambda',\mu'}$ is not the starting height of the sequence, we have $\lambda'$ starting with $i_1i_2\cdots i_pi_{p+1}i_{p+2}\cdots i_{p+r-1}$ and $\mu'$ starting with $i_2i_3 \cdots i_{p+1}i_{p+2} \cdots i_p$ for different indices $i_1,\ldots,i_{p+r-1}$. Then we can move the disjoint cycle $i_p \to i_{p+1} \to \ldots \to i_{p+r-1} \to i_p$ to the beginning, remove it, and repeat the algorithm with the rest of the sequence. The resulting sequence is a concatenation of disjoint cycles, and we call it the \emph{disjoint cycle decomposition} of the o-sequence $a_{\lambda,\mu}$ (or of the p-sequence $a_{\lambda',\mu'}$). For example, the disjoint cycle decomposition of
$$a_{13}a_{11}a_{12}a_{13}a_{22}a_{23}a_{22}a_{21}a_{23}a_{22}a_{23}a_{32}a_{31}a_{31}a_{33}a_{32}a_{32}a_{33}a_{33}$$
is
$$a_{22}a_{32}a_{23}a_{13}a_{31}a_{11}a_{22}a_{12}a_{21}a_{13}a_{31}a_{33}a_{23}a_{32}a_{22}a_{23}a_{32}a_{33}a_{33}.$$

\medskip

We say that two cycles in the disjoint cycle decomposition are \emph{disjoint} if the sets of their starting heights are disjoint.

\medskip

Recall that for a Cartier-Foata matrix $A$, the matrix $C=(c_{ij})_{n+1 \leq i,j \leq m}$ with
$$c_{ij} = - \det{}^{-1}(I-A_0) \cdot \det \begin{pmatrix} I-A_0 & -a_{*j} \\ -a_{i*} & -a_{ij} \end{pmatrix}$$
is right-quantum by Lemma \ref{cf2} and Proposition \ref{cf5}, so
$$\det{}^{-1}(I-C) = (I-C)^{-1}_{n+1,n+1}(I-C^{n+1,n+1})^{-1}_{n+2,n+2}\cdots = $$
\begin{equation} \label{beta4}
 = (I-A)^{-1}_{n+1,n+1}(I-A^{n+1,n+1})^{-1}_{n+2,n+2}\cdots
\end{equation}
by Theorem \ref{sylv2}. The last expression is the sum over all sequences which are concatenations of a lattice path from $n+1$ to $n+1$, a lattice path from $n+2$ to $n+2$, etc.

\begin{comment}
Since $A$ is Cartier-Foata, we can reorder this sequence so that the starting heights of steps are non-decreasing. As a running example, we will take $m=4$, $n=2$ and the o-sequence
$$a_{11}a_{13}a_{12}a_{25}a_{22}a_{21}a_{24}a_{32}a_{31}a_{43}a_{52}a_{55},$$
which is modulo $\p I_{\cf}$ equal to the concatenation of paths $3 \to 2 \to 5 \to 2 \to 2 \to 1 \to 1 \to 3 \to 1 \to 2 \to 4 \to 3$ and $5 \to 5$.

\medskip

Take an o-sequence $\alpha$ with the disjoint cycle decomposition $u_1u_2\ldots u_k$. Let $\p J$ be the set of all $i \in \set{1,\ldots,k}$ such that $u_i$ contains a height $> n$. Then $\alpha$ appears in $\det{}^{-1}(I-C)$ if and only if the following property is satisfied. For each $i = 1,\ldots,k$ there exists $j \geq i$ with $j \in \p J$, and if $i<j$ there exists $k$, $i < k \leq j$, such that $u_i$ and $u_k$ are not disjoint. For example, the disjoint cycle decomposition of the o-sequence above is
$$(a_{11})(a_{25}a_{52})(a_{22})(a_{13}a_{32}a_{21})(a_{12}
{1}, {2, 5}, {2}, {1, 3, 2}, {1, 2, 4, 3}, {5}}

\medskip

This statement has the following generalization.
\end{comment}

\begin{thm}[$\beta$-extension of Cartier-Foata Sylvester's identity] \label{beta3}
 Assume $A=(a_{ij})_{m \times m}$ is a Cartier-Foata matrix. For
 $$C=(c_{ij})_{n+1 \leq i,j \leq m} \qquad \mbox{with} \qquad c_{ij} = - \det{}^{-1}(I-A_0) \cdot \det \begin{pmatrix} I-A_0 & -a_{*j} \\ -a_{i*} & -a_{ij} \end{pmatrix}$$
 and for each $\beta \in \C$, the expression
 $$\left( \frac 1 {\det(I-C)} \right)^{\beta}$$
 is equal to
 $$\sum e_\mu(\beta) a_{\lambda,\mu},$$
 where $\mu$ runs over all words in the alphabet $\{1,\ldots,m\}$, $\lambda$ is the non-decreasing rearrangement of $\mu$, and $e_\mu(\beta)$ is a polynomial function of $\beta$ that is calculated as follows. Let $u_1u_2\ldots u_k$ be the disjoint cycle decomposition of $a_{\lambda,\mu}$. Let $\p J$ be the set of $i \in \set{1,\ldots,k}$ such that $u_i$ contains a height $> n$. Then
 \begin{equation} \label{beta2}
  e_\mu(\beta) = \sum_\pi \binom{\beta+l-1-d(\pi)}{l},
 \end{equation}
 where $d(\pi)$ denotes the number of descents of the subword of $\pi(1)\pi(2)\ldots\pi(k)$ composed of $\pi(i) \in \p J$, $l$ is the cardinality of $\p J$, and the sum is over all permutations $\pi \in S_k$ with the following properties:
 \begin{enumerate}
  \item \label{beta3a} if $i<j,\pi(i)>\pi(j)$, then $u_{\pi(i)},u_{\pi(j)}$ are disjoint;
  \item \label{beta3b} for each $i = 1,\ldots,k$ there exists $j \geq  i$ such that $u_{\pi(j)} \in \p J$;
  \item \label{beta3c} if $j \geq i$ is the minimal $j$ with $u_{\pi(j)} \in \p J$ and if $j > i$, there exists $k$, $i < k \leq j$, such that $u_{\pi(i)}$ and $u_{\pi(j)}$ are not disjoint;
  \item \label{beta3d} if $\pi(i)>\pi(i+1)$ then $\pi(i) \in \p J$.
 \end{enumerate}
\end{thm}

\begin{exm}
 Take $\mu = 132521421325$. The disjoint cycle decomposition of the o-sequence $a_{11}a_{13}a_{12}a_{25}a_{22}a_{21}a_{24}a_{32}a_{31}a_{43}a_{52}a_{55}$ is
 $$u_1u_2u_3u_4u_5u_6 = (a_{11})(a_{25}a_{52})(a_{22})(a_{13}a_{32}a_{21})(a_{12}a_{24}a_{43}a_{31})(a_{55}).$$
 We have $\p J = \set{2,4,5,6}$, the only permutations in $S_6$ that appear in the sum \eqref{beta2} are $213456, 213465, 261345$ with $d(213456)=0$, $d(213465)=1$, $d(261345)=1$. Therefore
 $$e_\mu(\beta) = \binom{\beta + 3}4 + 2 \binom{\beta+2}4 = \frac{\beta^4}{8}+\frac{5 \beta^3}{12}+\frac{3 \beta^2}{8}+\frac{\beta}{12}.$$
\end{exm}

\begin{exm}
 Take $n=0$. In this case $\p J = \set{1,\ldots,k}$, only the first condition is not vacuously true on $\pi$, and we get the $\beta$-extension of MacMahon Master Theorem, \cite[Theorem 10.5]{konvalinka}.
\end{exm}

It is clear that each term of $(\det(I-C))^{-\beta}$ is an o-sequence modulo $\p I_{\cf}$, and that the coefficients of o-sequences are polynomial functions of $\beta$. Therefore it is enough to prove the theorem for $\beta \in \N$, and this is an enumerative problem. We are given an o-sequence $a_{\lambda,\mu}$ and $\beta$ slots, and we have to calculate in how many ways we can choose terms of $(\det(I-C))^{-1}$ in each slot so that their product is, modulo $\p I_{\cf}$, equal to $a_{\lambda,\mu}$.

\medskip

We start the proof with a lemma.

\begin{lemma}
 All the steps in a cycle of the disjoint cycle decomposition must be placed in the same slot.
\end{lemma}
\begin{prf}
 This is proved in exactly the same way as the proof of \cite[Lemma 10.4]{konvalinka}, since all we used there was that the sequence chosen in each slot must be balanced, which is also true in our case. \qed
\end{prf}

\begin{proof}[Proof of Theorem \ref{beta3}] 
 We will call cycles with all heights $\leq n$ \emph{low cycles}, and cycles containing at least one height $>n$ \emph{high cycles}.\\
 The lemma tells us that we must choose a permutation $\pi \in S_k$ such that $u_1 \cdots u_k = u_{\pi(1)} \cdots u_{\pi(k)}$ modulo $\p I_{\cf}$, and place the cycles $u_{\pi(1)},\ldots,u_{\pi(k)}$ in the $\beta$ slots so that the cycles in each slot give a term appearing in $(\det(I-C))^{-1}$.\\
 Two cycles commute if and only if they are disjoint. Therefore the condition $u_1 \cdots u_k = u_{\pi(1)} \cdots u_{\pi(k)}$ implies \eqref{beta3a}.\\
 Take a low cycle $u_{\pi(i)}$, and assume that there are no high cycles $u_{\pi(j)}$ with $j > i$. That means that there are no high cycles to the right of $u_{\pi(i)}$ in its slot, and therefore the sequence in this slot is not of the form \eqref{beta4} modulo $\p I_{\cf}$, which is a contradiction. This proves \eqref{beta3b}.\\
 A permutation $\pi$ that satisfies these two conditions and not \eqref{beta3c}-\eqref{beta3d} may be allowable in the sense that there exists a placement of $u_{\pi(1)},u_{\pi(2)},\ldots,u_{\pi(k)}$ in this order in the $\beta$ slots so that the sequence in each slot is equal to a sequence of form \eqref{beta4} modulo $\p I_{\cf}$. However, different permutations can give exactly the same placements modulo $\p I_{\cf}$. For example, take the o-sequence $a_{11}a_{22}a_{31}a_{12}a_{23}$ for $m=3$, $n=2$. The disjoint cycle decomposition is $u_1u_2u_3 = (a_{11})(a_{22})(a_{31}a_{12}a_{23})$. Both $123$ and $213$ satisfy \eqref{beta3a}--\eqref{beta3b} in Theorem \ref{beta3}, but since $u_1,u_2,u_3$ must necessarily be in the same slot, these two permutations give exactly the same placements.\\
 We will show that all placements corresponding to permutations satisfying \eqref{beta3a}--\eqref{beta3b} also correspond to permutations satisfying \eqref{beta3a}--\eqref{beta3d}, and then count the number of placements corresponding to such permutations.\\
 Take a low cycle $u_{\pi(i)}$, and take the smallest $j>i$ so that $u_{\pi(j)}$ is a high cycle. If $u_{\pi(i)}$ is disjoint with $u_{\pi(k)}$ for $i < k \leq j$, we can move it to the right of $u_{\pi(j)}$, and hence reduce the number of violations of \eqref{beta3c}. Therefore we can assume that $\pi$ satisfies \eqref{beta3a}--\eqref{beta3c}.\\
 Finally, assume that we have $\pi(i)>\pi(i+1)$ for a low cycle $u_{\pi(i)}$. Then $u_{\pi(i+1)}$ must be placed in the same slot as $u_{\pi(i)}$ (otherwise the sequeunce in the slot containing $u_{\pi(i)}$ would not be equal to a sequence of the form \eqref{beta4} modulo $\p I_{\cf}$), and $u_{\pi(i)}$ and $u_{\pi(i+1)}$ commute, so we can switch $\pi(i)$ and $\pi(i+1)$ and reduce the number of violations of \eqref{beta3d}. Therefore we can assume that $\pi$ satisfies \eqref{beta3a}--\eqref{beta3d}.\\
 Now we have to find the number of ways to place $u_{\pi(1)},\ldots,u_{\pi(k)}$ in the $\beta$ slots so that the cycles in each slot give a term appearing in $(\det(I-C))^{-1}$. All cycles between two consecutive high cycles must appear in the same slot as the right-hand high cycle. Therefore placing the cycles in slots is the same as placing $\beta-1$ dividers after (some of the) high cycles. Of course, there are $\binom{\beta-1+l}l$ ways of doing this, but we can get the same terms several times: if we take two consecutive high cycles $u_{\pi(i)},u_{\pi(j)}$ with $i<j,\pi(i)>\pi(j)$, then $u_{\pi(i)}$ must necessairly commute with $u_{\pi(j)}$ \emph{and} with all the low cycles between them, we can move $u_{\pi(j)}$ to the right of $u_{\pi(i)}$, possibly move some of the low cycles before $u_{\pi(j)}$ to the right of $u_{\pi(i)}$, and we see that this term has already been counted by a different $\pi$. In order to avoid overcounting, we \emph{have to} place a divider after $u_{\pi(i)}$. Therefore the number of unique placements in slots corresponding to $\pi$ is $\binom{\beta-1+l-d(\pi)}l$, and this finishes the proof of Theorem \ref{beta3}.
\end{proof}

\section{Final remarks} \label{final}

\subsection{}
The paper \cite{konvalinka} gives proofs of only special cases of Propositions \ref{cf1}, \ref{cfq1} and \ref{cfqij9}. However, the proof given in \cite[Section 12]{konvalinka} can be easily extended. Namely, for $B=I-A$ with $A$ Cartier-Foata (resp.\hspace{-0.07cm} right-quantum), the $j$-th coordinate of the matrix product
$$((-1)^{i+1}\det B^{1i},(-1)^{i+2} \det B^{2i},\ldots,(-1)^{i+m} B^{mi}) \cdot B$$
is $\sum_{k=1}^m (-1)^{i+k} \det B^{ki} b_{kj}$; since $B$ satisfies the conditions of \cite[Lemma 12.1]{konvalinka} (resp.\hspace{-0.07cm} \cite[Lemma 12.2]{konvalinka}), this is equal to $\det B \cdot \delta_{ij}$. Then
$$((-1)^{i+1}\det B^{1i},(-1)^{i+2} \det B^{2i},\ldots,(-1)^{i+m} B^{mi}) = \det B \cdot (0,\ldots,1,\ldots,0) \cdot B^{-1}$$
and
$$\left(B^{-1}\right)_{ij} = (-1)^{i+j} \det{}^{-1} B \cdot B^{ji}.$$

Furthermore, if $i = m$, we never have to use any conditions on elements of the form $a_{km}$: $\sum_{k=1}^m (-1)^{m+k} \det B^{km} b_{km} = \det B$ by part (3) of \cite[Lemma 12.1]{konvalinka} (resp.\hspace{-0.07cm} \cite[Lemma 12.2]{konvalinka}) and no element of $\sum_{k=1}^m (-1)^{m+k} \det B^{km} b_{kj}$ for $j < m$ is of the form $b_{km}$. 

\medskip

This proves Propositions \ref{cf1}, and the $1=q$ and $1=q_{ij}$ principles (\cite[Lemma 12.3]{konvalinka} and \cite[Lemma 12.4]{konvalinka}) prove Propositions \ref{cfq1}, \ref{cfqij9} and \ref{cfqij10}. See also \cite{konvalinka2} for a combinatorial proof.

\subsection{}
Bareiss's proof of Theorem \ref{intro2} is a pretty straighforward linear algebra argument; see \cite{muhlbach}, \cite{akritas} for other proofs and some mild generalizations.

\medskip

In 1991, a generalization to quasideterminants, essentially equivalent to our Theorem \ref{sylv2}, was found by Gelfand and Retakh \cite{gelfand1}. Krob and Leclerc \cite{krob} used their result to prove the following quantum version.

\medskip

Let $q \in \C \setminus \set 0$. Call a matrix (in non-commutative variables) $A=(a_{ij})_{m \times m}$ quantum if:
\begin{itemize}
 \item $a_{jk}a_{ik} = q a_{ik} a_{jk}$ for $i < j$,
 \item $a_{il}a_{ik} = q a_{ik} a_{il}$ for $k < l$,
 \item $a_{jk}a_{il} = a_{il}a_{jk}$ for $i < j,k < l$,
 \item $a_{ik}a_{jl} - a_{jl}a_{ik} = (q^{-1} - q) a_{il}a_{jk}$ for $i < j,k < l$.
\end{itemize}
Define the quantum determinant of a matrix $B$ by
$$\detq A = \sum_{\sigma \in S_m} (-q)^{-\inv \sigma} a_{\sigma(1)1}a_{\sigma(2)2}\cdots a_{\sigma(m)m},$$
where $\inv \sigma$ denotes the number of inversions of the permutation $\sigma$.

\begin{thm}[Krob, Leclerc]
 For a quantum matrix $A=(a_{ij})_{m \times m}$, take $n$, $A_0$, $a_{i*}$ and $a_{*j}$ as before, and define
 $$b_{ij} = \detq \begin{pmatrix} A_0 & a_{*j} \\ a_{i*} & a_{ij} \end{pmatrix}, \quad B = (b_{ij})_{n+1 \leq i,j \leq m}.$$
 Then
 $$\detq A \cdot (\detq A_0)^{m-n-1} = \detq B. \eqno \qed$$
\end{thm}

Krob and Leclerc's proof consists of an application of the so-called quantum Muir's law of extensible minors to the expansion of a minor.

\bigskip

{\bf Acknowledgements.} The author is grateful to Igor Pak for numerous suggestions, and to Pavel Etingof and Alexander Molev for help with the references.

\medskip

{\sc \scriptsize Department of Mathematics, Massachusetts Institute of Technology, Cambridge, MA 02139\\
\tt{konvalinka@math.mit.edu}\\
\tt{http://www-math.mit.edu/\~{}konvalinka/}}
                                                                                           
\end{document}